\documentclass[conference]{IEEEtran}

\IEEEoverridecommandlockouts                              

\usepackage{graphicx} 
\usepackage{amssymb}
\usepackage{amsmath}
\usepackage{mathrsfs}
\usepackage{algorithm}
\usepackage{algorithmic}
\usepackage{amsthm}
\usepackage{xcolor}

\graphicspath{ {./images/} }
\newtheorem{theorem}{Theorem}
\newtheorem{lemma}{Lemma}
\newtheorem{definition}{Definition}
\newtheorem{assumption}{Assumption}
\newtheorem{remark}{Remark}

\makeatletter
\newenvironment{breakablealgorithm}
  {
    \begin{center}
      \refstepcounter{algorithm}
      \hrule height.8pt depth0pt \kern2pt
      \parskip 0pt
      \renewcommand{\caption}[2][\relax]{
        {\raggedright\textbf{\fname@algorithm~\thealgorithm} ##2\par}%
        \ifx\relax##1\relax 
          \addcontentsline{loa}{algorithm}{\protect\numberline{\thealgorithm}##2}%
        \else 
          \addcontentsline{loa}{algorithm}{\protect\numberline{\thealgorithm}##1}%
        \fi
        \kern2pt\hrule\kern2pt
     }
  }
  {
     \kern2pt\hrule\relax
   \end{center}
  }
\makeatother

\title{Distributed Least-Squares Optimization Solvers with Differential Privacy}
\author{Weijia Liu, Lei Wang, Fanghong Guo, Zhengguang Wu, Hongye Su
 \thanks{This work  was supported in part by the National
Key R\&D Program of China under Grant 2018YFA0703800, in
part by Zhejiang Provincial Natural Science Foundation of China under Grant No. LZ23F030008 and in part by the National Natural Science Foundation of China under Grant No. 62203386 and 62373328.  }
\thanks{W. Liu, L. Wang, Z. Wu and H. Su are with  Institute of Cyber-Systems and Control, College of Control Science and Engineering, Zhejiang University, Hangzhou, 310027, China (Emails: liuweijia0430@163.com; lei.wangzju@zju.edu.cn; nashwzhg@zju.edu.cn; hysu@iipc.zju.edu.cn). F. Guo is with Department of Automation, Zhejiang University of Technology, Hangzhou, China (Emails: fhguo@zjut.edu.cn).}
}

\begin{document}

\maketitle
\thispagestyle{empty}
\pagestyle{empty}
\begin{abstract}
This paper studies the distributed least-squares optimization problem with differential privacy requirement of local cost functions, for which two differentially private distributed solvers are proposed. 
The first is established  on the distributed gradient tracking algorithm, by appropriately perturbing the initial values and parameters that contain the privacy-sensitive data with Gaussian and truncated Laplacian noises, respectively. Rigorous proofs are established to show the achievable trade-off between the $(\epsilon,\delta)$-differential privacy and the computation accuracy.  
The second solver is established on the combination of the distributed shuffling mechanism and the average consensus algorithm, which enables each agent to obtain a noisy version of parameters characterizing the global gradient. As a result,  the least-squares optimization problem can be eventually solved by each agent locally in such a way that any given $(\epsilon,\delta)$-differential privacy requirement can be preserved while the solution may be computed with the accuracy independent of the network size, which makes the latter more suitable for large-scale distributed least-squares problems. 
Numerical simulations are presented to show the effectiveness of both solvers.
\end{abstract}

\section{Introduction}
Distributed optimization, in which networked agents seek to minimize a global cost function in a distributed and cooperative fashion, has gained significant attention in diverse fields such as machine learning \cite{DisOptML}, unmanned aerial vehicles (UAV) \cite{DisOptUAV}, and sensor networks \cite{DisOptSensor}. For such a problem, since  the seminal work of the distributed gradient-descend algorithm \cite{DisOptGD}, numerical  algorithms have been developed, that differ from linear or sub-linear convergence  rate \cite{GradientTracking, DisOptNesterovGradient,DisOptSubgradientPush,DisOptPI},  deterministic or stochastic communication graph \cite{DisOptStochasticGradient,DisOptAMMD}, and equality or inequality constraints \cite{DisOptProjectedSubgradient,DisOptLagrangian}, etc. Particularly, the gradient-tracking algorithm \cite{varagnolo2015newton,GradientTracking} has gained increasing attention due to its fast convergence rate, yielding several variants 
 such as  \cite{nedic2017achieving,xu2017convergence,scutari2019distributed,pu2020push,pu2021distributed}. As a special case, distribution least-squares optimization considers the computation problems where  cost functions are squared, with wide applications across many  fields (e.g., versatile LiDAR SLAM \cite{LSSLAM}). 
 See \cite{DisOptReview,yang2019survey} for a thorough review of recent advances of distributed optimization.  


In distributed optimization, each network agent holds a local dataset (that defines a local cost function), that  may  contain private sensitive information for the agent. For example, in the digital contact tracing based computational epidemiology, users' privacy-sensitive data such as localization is encoded in the corresponding optimization problem \cite{DisOptPrivacyLeakage}. In existing distributed optimization algorithms, during computing processes each agent needs to communicate with neighboring agents its local states or some intermediate variables that may directly contain the local sensitive data, or can be used to infer the local data.  In view of this, new risks of privacy breach arise when there are eavesdroppers having access to network communication messages.

For privacy concern in distributed optimization, a natural idea is to apply homomorphic encryption on transmitted messages \cite{DisOptAMMD, DisOptHorEncryp1}. However, this method may lead to heavy computation and communication burden. Besides, other approaches have also been proposed by injecting uncertainties into exchanged states \cite{DisOptPrivacyNoiseStates1, YQWangDisOptDP}, or into asynchronous heterogeneous stepsize \cite{DisOptPrivacyProjStep}. To further describe a rigorous privacy preserving notion, the differential privacy is proposed in \cite{DworkDP}, and has been widely applied in various fields such as signal processing \cite{DPSignalProcessing2},  and distributed computation \cite{YQWangDisOptDP,LapDPAC,Dishuf,DPDisOptMess1,DPDisOptMess2,DPDisOptFunc,DPDisOptHan2017,DPDisOptCloudHale2018}, due to its resilience to side information and post-processing \cite{DPpost-processing}. 
For a distributed optimization problem with a trusted cloud,  a differentially private computation framework is proposed in \cite{DPDisOptCloudHale2018} where the cloud sends perturbed global gradient to agents. 
Further, for the scenario without a trusted cloud, \cite{DPDisOptMess1} studies  distributed optimization problems with bounded gradients and sensitive constraints by  injecting noises to communication messages, for which the desired privacy is guaranteed by splitting the privacy budget to each iteration round.   In \cite{DPDisOptFunc}, to preserve the differential privacy under arbitrary iteration rounds, the idea of perturbing objective functions is introduced for a class of distributed convex optimization problem with square-integrable cost functions. 
Note that as shown  in \cite{DPDisOptFunc}, in these results the introduction of randomness for privacy concern inevitably leads to a certain computation error, i.e., there is an accuracy-privacy trade-off. When an exact optimal solution is sought, \cite{YQWangDisOptDP,DPDisOptMess2} turn to make a bit modification on the adjacency notion for the differential privacy, and develop novel algorithms by perturbing communication messages with noises having increasing variants and employing a vanishing stepsize that is appropriately designed.

In this paper, we focus on the problem of distributed least-squares optimization with differential privacy requirements of local cost functions under a common data adjacency notion. For such a problem, a common approach is to perturb communication messages, but as in \cite{DPDisOptMess1} the change of iteration number may lead to a different privacy budget, i.e., yielding an iteration-number-dependent differential privacy.  Besides, we also note that the idea of functional perturbation \cite{DPDisOptFunc} is not applicable as the cost functions are not square-integrable, though it may remove the dependence on the iteration number. 

In view of this analysis, inspired by the truncated Laplace mechanism \cite{BoundedNoiseDP} and the Gaussian mechanism \cite{ImprovedGaussianDP}, we first propose a Differentially Private Gradient Tracking-based (DP-GT-based) solver by perturbing the gradient tracking algorithm in such a way of appropriately adding Gaussian and truncated Laplacian noises to the initial values and parameters that contain the sensitive data, respectively. Moreover, a trade-off between 
the achievable $(\epsilon,\delta)$-differential privacy and the computation accuracy is rigorously proved, and  the iteration round is allowed to be arbitrary with no influence on the privacy guarantee. For such a method, it is noted that the intrinsic property of truncated Laplacian differential privacy mechanism \cite{BoundedNoiseDP} limits the achievable differential-privacy level and the computation accuracy linearly relies on the network size, which makes it inapplicable to a large-scaled distributed least-squares problem.
In view of this, we further propose a  Differentially Private Dishuf Average Consensus-Based (DP-DiShuf-AC-based) solver, where agents run the Distributed Shuffling (DiShuf)-based differentially private average consensus algorithm \cite{wang23DishufImprovedDPAC} to obtain a noisy estimate
of the global gradient function, enabling  each agent  to solve for the optimal solution independently with differential privacy guarantees. 
Benefiting from the use of distributed shuffling mechanism, the least-squares optimization problem can be eventually solved to fulfill any desired differential privacy levels, and with the accuracy independent of the network size, which makes it more suitable for large-scale least-squares optimization problems. 
In the paper, our contribution mainly lies in proposing two differentially private least-squares optimization solvers: the DP-GT-based solver and DP-DiShuf-AC-based Solver. Particularly, the latter is able to achieve a privacy-accuracy trade-off which is independent of the network size as shown by numerical simulations. 

The remainder of this paper is organized as follows. In Section \ref{sec-2}  the differentially private distributed least squares problem is formulated, for which  two solvers are explicitly presented in Sections \ref{sec-3} and \ref{sec-4}. Simulations are then conducted in Section \ref{sec-5} to verify the effectiveness of the proposed solvers. Finally, Section \ref{sec-6} concludes the paper.

{\bf Notation}. Denote by $\mathbb{R}$ the real number, $\mathbb{R}^{n}$ the real space of $n$ dimension for any positive integer $n$. 
Denote by $\mathbf{1}_{n}\in\mathbb{R}^{n}$ as a vector with each entry being $1$ and $\mathbf{0}_{n}\in\mathbb{R}^{n}$ as a vector with each entry being $0$.
For a vector $x\in\mathbb{R}^{n}$, denote ${\|x\|}_0$, ${\|x\|}_1$, ${\|x\|}$ as the 0, 1, and 2-norm of vector $x$. 
For a matrix $A$, denote $(A)_{ij}$ as its $i$-th row $j$-th column element, $\|A\|$ as 2-norm, $\|A\|_{F}$ as Frobenius norm, and $\lambda_i(A)$ as $i$-th eigenvalue.
Denote by $\eta\sim\mathcal{N}(\mu,\sigma^2)^r$ if each entry of $\eta\in\mathbb{R}^r$ is i.i.d. and drawn from Gaussian distribution with mean  $\mu$ and variance $\sigma^2$ for positive integer $r$. 
Define $\Phi(s)=\frac{1}{\sqrt{2\pi}}\int_{-\infty}^{s}e^{-\tau^2/2}d\tau$ and $\kappa_\epsilon(s)=\Phi(\frac{s}{2}-\frac{\epsilon}{s})-e^{\epsilon}\Phi(-\frac{s}{2}-\frac{\epsilon}{s})$ for any $\epsilon>0$. Denote by $\bar\kappa_\epsilon(\delta)$ the inverse function of $\kappa_\epsilon(s)$ for $s>0$. It can be verified that both $\kappa_\epsilon(s)$ and $\bar\kappa_\epsilon(\delta)$ are strictly increasing functions.

\section{Problem Statement}
\label{sec-2}
Consider a multi-agent system with an undirected and connected communication network $\mathrm{G}=(\mathrm{V},\mathrm{E})$, where the agent set $\mathrm{V}=\{1,2,...,n\}$, the edge set $\mathrm{E}\subseteq \mathrm{V}\times \mathrm{V}$, and each agent $i\in\mathrm{V}$ holds a local and \emph{privacy-sensitive} cost function of the form
\begin{equation}\label{eq:f_i}
    f_i(x)=\frac{1}{2}x^{\top} A_{i}x+B_{i}^{\top} x+C_{i},
\end{equation}
where $A_{i}\in \mathbb{R}^{m \times m}$ is symmetric, $B_{i}\in \mathbb{R}^{m}$ and $C_{i}$ is a scalar. 
The networked system aims to solve the  following least-squares optimization problem
\begin{equation}\label{eq:LS}
    \min_{x\in \mathbb{R}^{m}} f(x)=\sum_{i=1}^{n}f_i(x)
\end{equation}
in a \emph{distributed} and \emph{privacy-preserved} manner. It is well-known that the value of $C_i$ does not affect the solution of the optimization problem and is generally not used. Thus, in the sequel  we only consider the privacy concern of the matrix pair $(A_i,B_i)$ while developing algorithms to solve the problem \eqref{eq:LS}, though the privacy of $f_i$ is encoded in the triplet $(A_i,B_i,C_i)$. Denote by $\mathcal{F}$ as the mapping that rearranges all elements of matrix $A_i$ to a vector $\theta_{i}^A\in\mathbb{R}^{\frac{m(m+1)}{2}}$, i.e. $\theta_{i}^A=\mathcal{F}(A_i)$. Specifically, $\theta_{i}^A$ is obtained by setting the element in the $p$-th row and the $q$-th column with $q\geq p$ of $A_{i}$ as the $[\frac{(2m-p+2)(p-1)}{2}+q-(p-1)]$-th element.  Note that since $A_i$ is symmetric, we only put the  upper-triangle elements of $A_{i}$ in the mapping $\mathcal{F}$. On the other hand, it is clear that $\mathcal{F}$ is invertible. For convenience, we  denote by $\mathcal{F}^{-1}$ as its inverse, i.e., $A_{i}=\mathcal{F}^{-1}(\theta_{i}^A)$, and define the sensitive-data vector 
$\theta_i:=[{\theta_i^A}; B_i]\in\mathbb{R}^{\frac{m(m+3)}{2}}$. 

To ensure the solvability of the problem \eqref{eq:LS}, we make the following assumption throughout the paper.
\begin{assumption}\label{ass-1}
    The matrix $A:=\sum_{i=1}^{n}A_{i}$ is positive definite, i.e., there exists a $\underline{\lambda}_{A}>0$ such that $A\geq\underline{\lambda}_{A}\mathbf{I}_m$.
\end{assumption}
The above assumption guarantees strict convexity of global cost function $f(x)$, admitting the unique least-squares solution
\begin{equation}\label{eq:x_ast}
    x^\ast=-A^{-1}B\,,
\end{equation}
where we define $B:=\sum_{i=1}^{n}B_{i}$.

In solving the distributed least-squares optimization problem \eqref{eq:LS} over the network $\mathrm{G}$, necessarily agents communicate with neighboring agents about local information, which may contain the sensitive data $(A_{i},B_{i})$ (i.e., vector $\theta_{i}$) directly or indirectly. If there exists adversaries having access to the communication messages, then the sensitive data $(A_{i},B_{i})$ may be inferred, leading to privacy leakage risks. 
This motivates to develop distributed algorithms that solve the problem \eqref{eq:LS} with privacy guarantees of the sensitive data $(A_{i},B_{i})$ against the adversaries eavesdropping all communication messages. 

For privacy concern, this paper focuses on the notion of differential privacy, originated from \cite{DworkDP} and is widely used in the field of distributed computation. Differential privacy describes rigorous privacy preserving in that eavesdroppers could not infer sensitive data by the way of ``differentiating". Any pair of data $\theta_{i},\theta_{i}'\in \mathbb{R}^{\frac{m(m+3)}{2}}$ is said to be $\mu$-adjacent with $\mu>0$, denoted by $(\theta_{i},\theta_{i}')\in Adj(\mu)$, if ${\|\theta_{i}-\theta_{i}'\|}_0=1$ and ${\|\theta_{i}-\theta_{i}'\|}_1\leq\mu$. We denote $\mathcal{M}$ as a mapping (or mechanism) from the sensitive data $\theta_{i}$ to the eavesdropped communication  messages, and give the following definition of  $(\epsilon, \delta)$-differential privacy \cite{DworkDP}.

\begin{definition}\label{def:DP}{\em (Differential Privacy)}.
Given $\epsilon\geq0, \delta\in(0,1)$, the distributed solver  preserves $(\epsilon, \delta)$-differential privacy of $\theta_{i}$ under $\mu$-adjacency, if for all $\mathbb{M}\subseteq range(\mathcal{M})$, there holds
\begin{equation}\label{eq:DP definition}
    \mathbb{P}(\mathcal{M}(\theta_{i})\in\mathbb{M})\leq e^{\epsilon}\mathbb{P}(\mathcal{M}(\theta_{i}')\in\mathbb{M})+\delta
\end{equation}
for any $(\theta_{i},\theta_{i}')\in Adj(\mu)$.
\end{definition}

Before the close of this section, we make the following claims on the communication network.
Denote by $\mathrm{N_{i}}=\{j|(i,j)\in\mathrm{E}\}$ as neighbor set of agent $i$, and  $w_{ij}$ the weight of edge $(i,j)$, satisfying $w_{ij}=w_{ji}\in(0,1)$ for $j\in \mathrm{N_{i}}$ and $w_{ji}=0$ for $j\notin \mathrm{N_{i}}$. Moreover, denote by $L$ the Laplacian matrix of the graph $\mathrm{G}$, satisfying $(L)_{ij}=-w_{ij},\enspace i\neq j$, and $(L)_{ii}=\sum_{k=1}^{n}w_{ik}\in(0,1),\enspace i\in \mathrm{V}$. 
Since the graph $\mathrm{G}$ is assumed to be connected, by \cite[Theorem 2.8]{mesbahi2010graph} and Gershgorin Circle
Theorem \cite{horn_johnson_1985} there holds $0= \lambda_1(L) < \lambda_2(L) \leq \ldots \leq \lambda_n(L) < 2$.

\section{Gradient-Tracking-based Solver with Differential Privacy}
\label{sec-3}

In this section, we modify the distributed gradient tracking algorithm \cite{GradientTracking} by injecting random noises to the sensitive data for the purpose of solving the distributed least-squares problem \eqref{eq:LS} while preserving differential privacy of $\theta_i$, $i\in\mathrm{V}$.

Let  $\bar\gamma=\frac{d}{\sqrt{n}m}\underline{\lambda}_{A}$ for $d\in(0,1)$. For any $\epsilon>0$,  we let 
\begin{equation}\label{eq:DP-GT slover delta}
  \mu=c \bar\gamma\,,\quad  1/2 > \delta \geq \frac{e^\epsilon -1}{2(e^{\epsilon/c} -1)}
\end{equation}
for $c\in(0,1)$, and choose $\sigma_\eta\geq {\mu}/{\bar\kappa_{\epsilon}(\delta)}$. 
With these parameters, we propose the Differentially Private Gradient Tracking-based (DP-GT-based) solver in Algorithm 1.

\begin{breakablealgorithm}
\caption{Differentially Private Gradient Tracking-based solver (DP-GT-based solver)}
\leftline{{\bf Input:} Data $\theta_{i}$, privacy budgets $\mu,\epsilon,\delta$, parameters $\bar\gamma, \sigma_\eta$.}
\begin{itemize}
    \item [1.] Each agent $i\in\mathrm{V}$ independently generates a vector of i.i.d. truncated Laplacian noises $\gamma_{i}\in\mathbb{R}^{\frac{m(m+1)}{2}}$, with each entry independently generated according to the  probability density function 
    \begin{equation}\label{eq:Truncated Laplacian Distribution}
        p(\gamma)=\left\{ 
        \begin{aligned}
        & \frac{\epsilon}{2\mu(1-e^{-{\epsilon\bar\gamma}/{\mu}})}e^{-{\epsilon|\gamma|}/{\mu}}, & \gamma\in[-\bar\gamma,\bar\gamma]
        \\ & 0, & otherwise
        \end{aligned}
        \right.,
    \end{equation}
    and i.i.d. Gaussian noises $\eta_{i}\sim\mathcal{N}(0,\sigma_\eta^2)^m$.
    \item [2.] Each agent $i\in\mathrm{V}$ computes $G_{i}=\mathcal{F}^{-1}(\theta_{i}^A+\gamma_{i})$ and $H_{i}=B_{i}+\eta_{i}$, and initializes the local states as
    \begin{equation}
        x_i(0)=  \mathbf{0}_{m} \,,\quad 
        s_{i}(0)=  H_i.
    \end{equation}
    \item [3.] For $t=0,1,2,\dots$, each agent $i\in\mathrm{V}$ sends  $(x_{i}(t),s_{i}(t))$ to neighboring agents $j\in\mathrm{N}_{i}$ and iterates the local states by following 
    \begin{equation}\label{eq:gt}
        \begin{aligned}
            x_{i}(t+1) = & x_i(t)+\sum_{j\in\mathrm{N}_{i}} w_{ij}(x_{j}(t)-x_i(t))-\beta s_{i}(t) \\
            s_{i}(t+1) = & s_i(t)+ \sum_{j\in\mathrm{N}_{i}} w_{ij}(s_{j}(t)-s_i(t))\\ & +G_i(x_i(t+1)-x_i(t)).
        \end{aligned}
    \end{equation}
\end{itemize}
\end{breakablealgorithm}

In Algorithm 1, two types of random noises (i.e., truncated Laplacian noises $\gamma_i$ and Gaussian noises $\eta_i$) are, respectively, added to the gradient tracking algorithm in the manner of directly perturbing the parameters $A_i$ (or $\theta_i^A$) and $B_i$, for the purpose of preserving $(\epsilon,\delta)$-differential privacy of these parameters. Particularly, the use of truncated Laplacian distribution \eqref{eq:Truncated Laplacian Distribution} is motivated by \cite{BoundedNoiseDP}. It can be easily verified that the noise generated according to the truncated Laplacian distribution \eqref{eq:Truncated Laplacian Distribution} has the mean zero and variance $\sigma_\gamma^2$ satisfying 
\begin{equation}
\begin{array}{rcl}
    \sigma_{\gamma}^2&=&\displaystyle\frac{\epsilon}{2\mu(1-e^{-{\epsilon\bar\gamma}/{\mu}})} \int_{-\bar\gamma}^{\bar\gamma} \gamma^2 e^{-\frac{\epsilon|\gamma|}{\mu}}d\gamma\\
    &=& \displaystyle\frac{1}{1-e^{-{\epsilon\bar\gamma}/{\mu}}}[2\frac{\mu^2}{\epsilon^2}-e^{-\frac{\epsilon\bar\gamma}{\mu}}(\bar\gamma^2+2\frac{\mu}{\epsilon}\bar\gamma+2\frac{\mu^2}{\epsilon^2})].
\end{array}
\end{equation}
Denote $\Omega_{A}=\mathcal{F}^{-1}(\sum_{i=1}^{n}\gamma_{i})$ and $\Omega_{B}=\sum_{i=1}^{n}\eta_{i}$. The following result demonstrates the differential privacy and computation accuracy that can be achieved by the above algorithm.

\begin{theorem}
Given any privacy budgets $(\epsilon,\delta,\mu)$ satisfying \eqref{eq:DP-GT slover delta}, by letting $\sigma_\eta\geq {\mu}/{\bar\kappa_{\epsilon}(\delta)}$, there exists a $\beta^\ast>0$ such that for all $\beta\in(0,\beta^\ast)$, Algorithm 1 achieves the following  properties.
    \begin{itemize}
        \item [i).] {\bf (Privacy)} The $(\epsilon,\delta)$-differential privacy of $\theta_i$,  $i\in\mathrm{V}$ is preserved under $\mu$ adjacency.
        \item [ii).] {\bf (Convergence)} There exists $\alpha_1\in(0,1)$ such that
 $\|x_{i}(t)-x(\infty)\|=\mathcal{O}(\alpha_1^t)$ for all $i\in\mathrm{V}$, with $x(\infty):=-(A+\Omega_{A})^{-1}(B+\Omega_{B})$.            
\item[iii).] {\bf (Accuracy)} The mean-square error between  $x(\infty)$ and $x^\ast$ satisfies
\[
    \mathbb{E}\|x(\infty)-x^{\ast}\|^2\leq \frac{2nm^2\sigma_\gamma^2\|x^\ast\|^2  + 2nm\sigma_\eta^2}{(1-d)^2\underline\lambda_{A}^2} \,.
\]

    \end{itemize}
\end{theorem}

\begin{remark}
\label{rem-1}
    The algorithm \eqref{eq:gt} indeed solves the least-squares problem $\min_{x\in \mathbb{R}^{m}} \bar f(x)=\sum_{i=1}^{n}\frac{1}{2}x^{\top}G_{i}x+H_{i}^{\top}x$. To preserve the convexity of $\bar f(x)$ while achieving the differential privacy, truncated Laplacian noises $\gamma_i$ are injected with the truncated level $\bar\gamma<\underline{\lambda}_{A}/{(\sqrt{n}m)}$. However, due to the presence of such a constraint on the truncated level, from \eqref{eq:DP-GT slover delta} it can be seen that the achievable privacy budgets $(\epsilon,\delta,\mu)$ cannot be arbitrarily prescribed. Moreover, one can see from the statement iii) that the mean-square error $\mathbb{E}\|x(\infty)-x^{\ast}\|^2$ is a linear function of the network size $n$, indicating a worse computation accuracy under a larger network size. All these issues thus motivate us to develop a new solver that can guarantee arbitrarily given privacy budgets $(\epsilon,\delta,\mu)$ while allowing a better computation accuracy when a large scale of networked least-squares problem is investigated.
\end{remark}

\section{Average-Consensus-Based Solver with Differential Privacy}
\label{sec-4}

In this section, a new distributed solver with differential privacy is proposed to overcome the issues addressed in Remark \ref{rem-1}.

We first observe that, the least-squares solution $x^\ast=-(\sum_{i=1}^{n}A_{i})^{-1}(\sum_{i=1}^{n}B_{i})$ implies that the least-squares problem \eqref{eq:LS} is solved if each agent obtains the global information $A=\sum_{i=1}^{n}A_{i}$ and $B=\sum_{i=1}^{n}B_{i}$. This intuition immediately motivates us to employ differentially private average consensus algorithms (e.g., \cite{LapDPAC,qiaoDishufDPAC,wang23DishufImprovedDPAC}) to compute these global information in a distributed manner, while preserving the desired differential privacy. Further, we note that the differentially private average consensus algorithm proposed in \cite{qiaoDishufDPAC,wang23DishufImprovedDPAC} yields a computation accuracy that can almost recover the centralized averaging algorithm where the accuracy of computing the sum is independent of the network size.

In view of the previous analysis, we propose a new solver by employing the  differentially private average consensus algorithm proposed in \cite{qiaoDishufDPAC,wang23DishufImprovedDPAC}, whose implementation relies on the following distributed shuffling (DiShuf) mechanism where Paillier cryptosystem \cite{PaillierCryptoSystem} is adopted with $\mathrm{E}_i(\cdot)$ and $\mathrm{D}_i(\cdot)$ as encryption and decryption operation, respectively, and $(k_{pi}$, $k_{si}$) as public and private key pair. 

\begin{breakablealgorithm}
\caption{Distributed Shuffling(DiShuf) Mechanism}
\leftline{{\bf Input:} Data $\theta_{i}\in\mathbb{R}^{{{m(m+3)}/{2}}}$, variance $\sigma_\eta^2$,  key pairs}
\leftline{\quad \quad \quad  ($k_{pi}$, $k_{si}$), a large positive integer $\bar a>>1$.}

    \begin{itemize}
        \item[1.] Each agent $i$ generates a vector of i.i.d. Gaussian noise $\eta_{i}\sim\mathcal{N}(0,\sigma_\eta^2)^{{m(m+3)}/{2}}$ and adds to $\theta_{i}$, i.e. $\bar \theta_{i}:=\theta_{i}+\eta_{i}$.
        \item[2.] Each agent $i$ encrypts each entry of $-\bar \theta_{i}$ with local public key $k_{pi}$, yielding a ciphertext vector $\mathrm{E}_i(-\bar \theta_{i})$, then sends $\mathrm{E}_i(-\bar \theta_{i})$ and public key $k_{pi}$ to neighboring agents $j\in\mathrm{N}_i$.
        \item[3.] Each agent $i$ encrypts each entry of $\bar \theta_{i}$ with received neighboring public keys $k_{pj}$, yielding a ciphertext vector $\mathrm{E}_j(\bar \theta_{i})$, and computes in entry to obtain a vector $c_{ij}=\mathrm{E}_j(\bar \theta_{i})\mathrm{E}_j(-\bar \theta_{j})$ for $j\in \mathrm{N}_i$.
        \item[4.] Each agent $i$ generate a positive integer $a_{i \rightarrow j}\in [\frac{\bar a}{\sqrt{2}},\bar a]$, and computes in entry to get $(c_{ij})^{a_{i \rightarrow j}}$ for each $j\in \mathrm{N}_{i}$.
        \item[5.] Each agent $i$ sends $(c_{ij})^{a_{i \rightarrow j}}$ to $j\in \mathrm{N}_i$ and decrypts received $(c_{ji})^{a_{j \rightarrow i}}$ with local private key $k_{si}$, yielding $\mathrm{D}_i((c_{ji})^{a_{j \rightarrow i}})$ for $j\in \mathrm{N}_i$.
        \item[6.] Each agent $i$ multiplies each $\mathrm{D}_i((c_{ji})^{a_{j \rightarrow i}})$ by $a_{i\rightarrow j}, j\in \mathrm{N}_i$ and computes the sum
        \begin{equation}
            \Delta_{i}=\sum_{j\in \mathrm{N}_i}a_{i\rightarrow j}\mathrm{D}_i((c_{ji})^{a_{j \rightarrow i}}).
        \end{equation}
    \end{itemize}
\leftline{{\bf Output:} $\Delta_{i},\enspace i\in\mathrm{V}$}
\end{breakablealgorithm}

With some lengthy but straightforward computations, it can be seen that the output of the DiShuf mechanism (see Algorithm 2) satisfies
\[
\Delta_{i} = \sum_{j\in\mathrm{N}_i}a_{i\rightarrow j}a_{j\rightarrow i}(\theta_j-\theta_i + \eta_j-\eta_i)\,,
\]
which implies $\sum_{i\in\mathrm{V}} \Delta_{i}=\textbf{0}_m$. Then we introduce the differentially private DiShuf-based average consensus algorithm below.

\begin{breakablealgorithm}
\caption{Differentially Private DiShuf Average Consensus-based solver (DP-DiShuf-AC-based solver)}
\leftline{{\bf Input:} Data $\theta_{i}\in\mathbb{R}^{{m(m+3)}/{2}}$, variance $\sigma_\gamma^2$, $\zeta>0$}
    \begin{itemize}
        \item[1.] Each agent $i$ runs the DiShuf mechanism with input data $\theta_i$ and outputs $\Delta_i$, then generates a vector of i.i.d. Gaussian noises $\gamma_{i}\sim\mathcal{N}(0,\sigma_\gamma^2)^{{m(m+3)}/{2}}$, and initializes its state as
        \begin{equation}
            y_{i}(0)=\theta_{i}+\zeta\Delta_{i}+\gamma_{i}.
        \end{equation}
        \item[2.] For $t=0,1,2,\dots$,  each agent $i$ sends its state $y_i(t)$ to its neighboring agents $j\in\mathrm{N}_{i}$ and updates its state following
                \begin{equation}\label{eq:state iteration}
                    y_{i}(t+1)=y_{i}(t)+\sum_{j\in \mathrm{V}}w_{ij}(y_{j}(t)-y_{i}(t)).
                \end{equation}
         \end{itemize}
\end{breakablealgorithm}
As $\sum_{i=1}^n\Delta_{i}=\textbf{0}_m$, it can be easily verified  that 
\[
\sum_{i=1}^ny_i(0) = \sum_{i=1}^n \theta_i + \sum_{i=1}^n \gamma_i\,,
\]
implying that the noises $\eta_i$ do not affect the average of $y_i(0)$, $i\in\mathrm{V}$. Thus, there are two design freedoms of noise variances (i.e., $\sigma_\eta^2$ and $\sigma_\gamma^2$), with the former affecting only achievable differential privacy, but having no influence on the average.
As a result,  by appropriately designing $\sigma_\eta$ and $\sigma_\gamma$,  a better trade-off between the privacy and the accuracy can be achieved. The following result from \cite{wang23DishufImprovedDPAC} demonstrates the corresponding privacy and computation properties.

\begin{lemma}
\cite[Theorem 3]{wang23DishufImprovedDPAC}
Suppose Assumption \ref{ass-1} holds and let $\zeta=\frac{1}{n\bar a^2 +1}$ and $g>0$. Given any privacy budgets $\epsilon>0,\delta\in(0,1),\mu>0$, 
let
\begin{equation}
    \begin{array}{rcl}
         \sigma_\gamma &=& \frac{(1+g)\mu}{\sqrt{n}\bar\kappa_{\epsilon}(\delta)}\\
         \sigma_\eta &=& \frac{(n-1)\alpha^2}{(1-\alpha)^2(\bar\kappa_{\epsilon}(\delta))^2}[\frac{(1+g)^2\mu^2}{(1+g)^2-1}-\frac{(1+g)^2\mu^2}{n(n-1)\alpha^2}] 
    \end{array}
\end{equation}
 where $\alpha=(1-\frac{1}{(2(n+\bar a^{-2}))^{n-1}})^{1/(n-1)}$,  Algorithm 3 achieves the following  properties.
    \begin{itemize}
        \item [i).] {\bf (Privacy)} The $(\epsilon,\delta)$-differential privacy of $\theta_i$,  $i\in\mathrm{V}$ is preserved under $\mu$ adjacency.
        \item [ii).] {\bf (Convergence)} There holds
 $\|y_{i}(t)-y(\infty)\|=\mathcal{O}(\alpha_2^t)$ for all $i\in\mathrm{V}$, with $y(\infty):=\sum_{i=1}^n \theta_i/n + \sum_{i=1}^n \gamma_i/n$ and $\alpha_2=\max\{|1-\lambda_2(L)|,|1-\lambda_n(L)|\} < 1$.         
\item[iii).] {\bf (Accuracy)} The mean-square error between  $y(\infty)$ and $y^\ast:=\sum_{i=1}^n \theta_i/n$ satisfies
\[
    \mathbb{E}\|y(\infty)-y^{\ast}\|^2 = \frac{(1+g)^2 \mu^2}{n^2(\bar\kappa_\epsilon(\delta))^2} \,.
\]

    \end{itemize}
\end{lemma}

Towards this end, by employing Algorithm 3 to compute the average of  
$\theta_i:=[{\theta_i^A}; B_i]$, $i\in\mathrm{V}$ in a distributed and differentially private manner, each agent $i\in\mathrm{V}$ can eventually obtain $\hat\theta: = ny(\infty):=\sum_{i=1}^n \theta_i + \sum_{i=1}^n \gamma_i$ and thus an estimate of the global information $(A,B)$ as $\hat A = A+\Omega_A$ and $\hat B = B + \Omega_B$, with $\Omega_A=\mathcal{F}^{-1}(\omega_A)$ where $\omega_A$ denotes the vector formed by the first $m(m+1)/2$ entries of $\sum_{i=1}^n \gamma_i$, and $\Omega_B$ the vector formed by the last $m$ entries of $\sum_{i=1}^n \gamma_i$. Therefore, each agent $i$ can compute for the solution $x^\ast$ of the least-squares problem \eqref{eq:LS} by computing $\hat x^\ast$ from
    \begin{equation}
       \hat A \hat x^\ast=-\hat B.
    \end{equation}
It is noted that with $A$ invertible, its noisy version $\hat A = A+\Omega_A$ is actually invertible in a generic sense. Besides, if the resulting $\hat A$ is singular, one can implement Algorithm 3 again.
It is also noted that $\mathbb{E}\|\hat\theta - \sum_{i=1}^n \theta_i\|^2 = \frac{(1+g)^2 \mu^2}{(\bar\kappa_\epsilon(\delta))^2}$, which is independent of the network size $n$. A natural conjecture  comes out that the mean-square error between the computed $\hat x^\ast$ and the actual solution $x^\ast$ is also independent of $n$, though it is still open to derive an explicit expression of $\mathbb{E}\|\hat x^\ast - x^\ast\|^2$ due to the inverse operation on the random matrix, i.e., $\hat A^{-1}$.

\section{Case Studies}
\label{sec-5}
In this section, we implement numerical simulations to show the effectiveness of our proposed Algorithms 1 and 3 (i.e., DP-GT-based solver and DP-DiShuf-AC-based solver, respectively), and compare their computation performance. We  consider a cycle communication graph with each edge having the same weight $w_{ij}=0.3$ for $(i,j)\in\mathrm{E}$. In this section, $A_i\in\mathbb{R}^{3\times 3}$ and $B_i\in\mathbb{R}^3$ are randomly generated and  $A=\sum_{i=1}^{n}A_{i}>0$.

\begin{figure}
    \centering
    \includegraphics[width=0.45\textwidth]{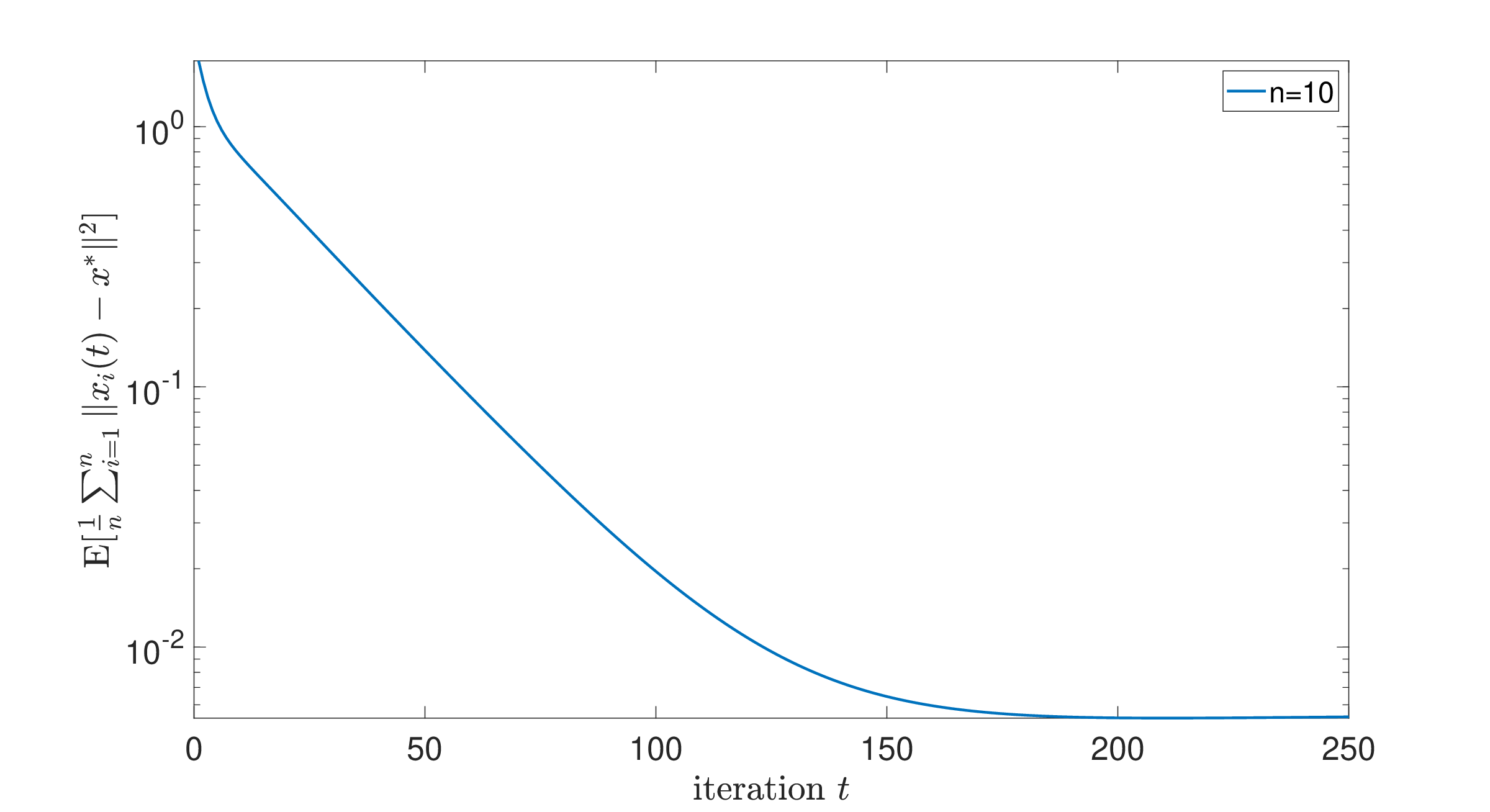}
    \caption{Trajectories of mean-square computation errors $\mathbb{E}[\frac{1}{n}\sum_{i=1}^{n}\|x_i(t)-x^\ast\|^2]$ of Algorithm 1}
    \label{fig:DP-GT error-t}
\end{figure}

Let the network size $n=10$ and choose $\bar\gamma=3.1$. Given the privacy budget $\mu=3$, $\epsilon=10$ and $\delta=0.2$, we let $\sigma_{\eta} = \frac{\mu}{\bar\kappa_{\epsilon}(\delta)}$ and $\beta=0.005$. It is shown in Figure \ref{fig:DP-GT error-t} that the mean-square error $\mathbb{E}[\frac{1}{n}\sum_{i=1}^{n}\|x_{i}(t)-x^{\ast}\|^{2}]$ decreases and exponentially converges as $t$ increases, demonstrating the effectiveness of Algorithm 1.

\begin{figure}
    \centering
    \includegraphics[width=0.45\textwidth]{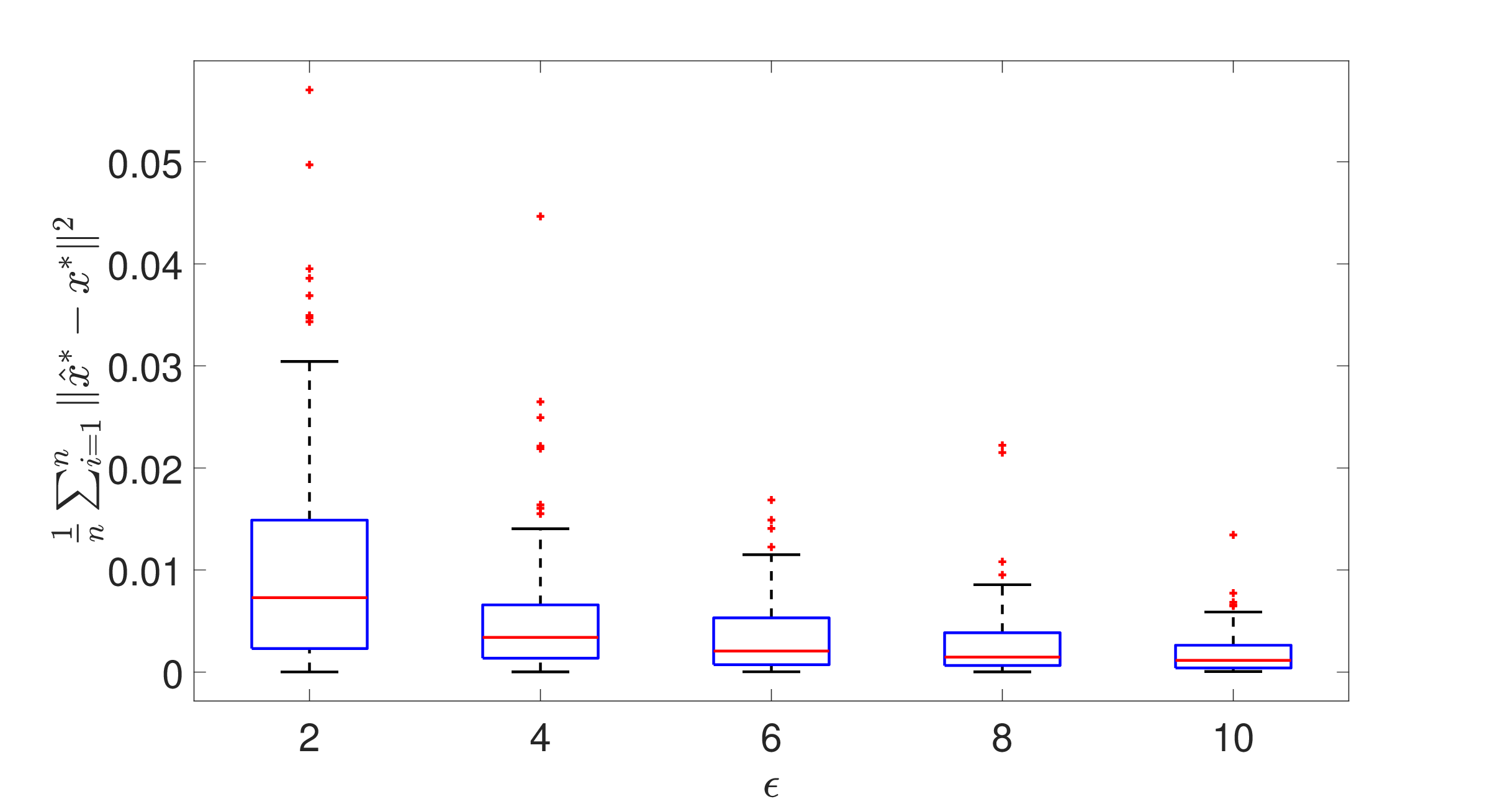}
    \caption{Box chart of computation errors $\frac{1}{n}\sum_{i=1}^n\| \hat x^\ast - x^\ast \| ^ 2$ of Algorithm 3 under different $\epsilon$ with 100 samples}
    \label{fig:Dishuf AC error-epsilon}
\end{figure}

Then we show the effectiveness of DP-DiShuf-AC-based solver. Taking the same $m$, $n$, $\mu$, $\delta$, choose $g=0.01$, we change $\epsilon$ and corresponding computed $\sigma_\gamma$ and $\sigma_\eta$. Figure \ref{fig:Dishuf AC error-epsilon} shows the distribution of computation error $\frac{1}{n}\sum_{i=1}^n\| \hat x^\ast - x^\ast \| ^ 2$ under different $\epsilon$ with 100 samples. It is observed our Algorithm 3 holds a resilience to a higher privacy preservation level, i.e., a smaller $\epsilon$, in the sense of computation accuracy.

\begin{figure}
    \centering
    \includegraphics[width=0.48\textwidth]{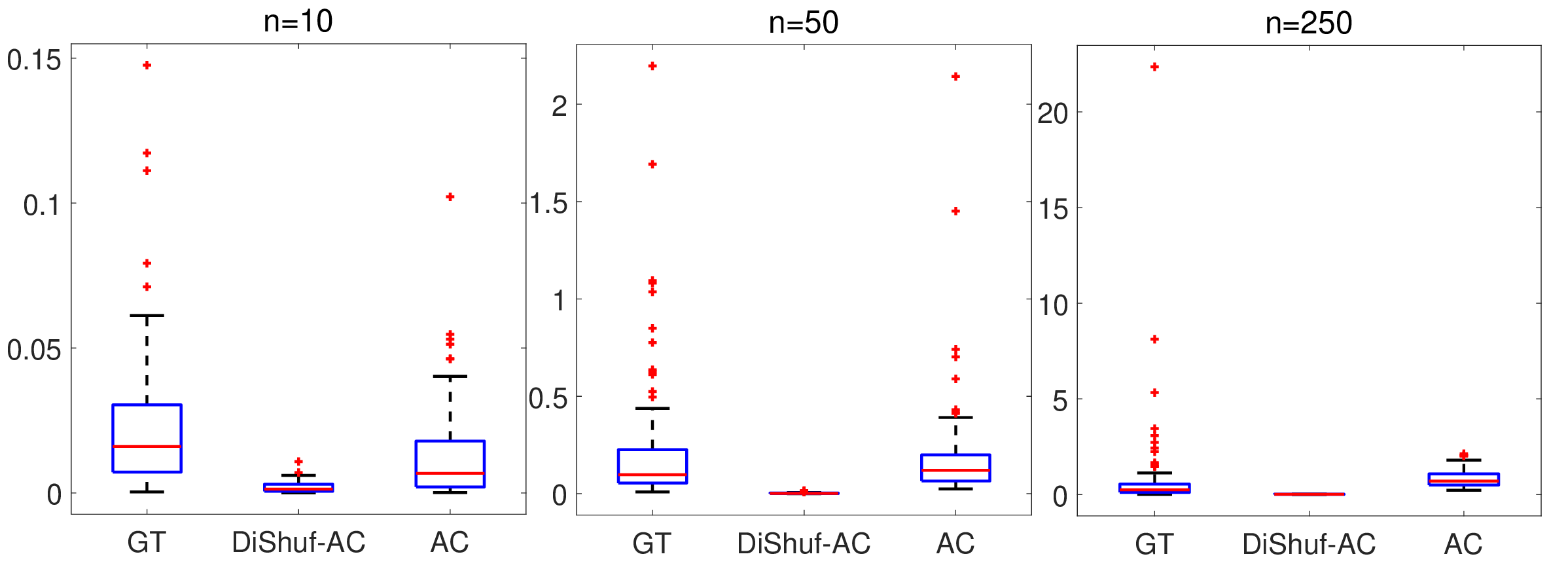}
    \caption{Box chart of computation errors $\frac{1}{n}\sum_{i=1}^{n}\| \hat x^\ast - x^\ast \|^2$ of three algorithms under $n=10,50,250$ with 100 samples}
    \label{fig:n comparison}
\end{figure}

Further, we study the computation error $\frac{1}{n}\sum_{i=1}^{n}\| \hat x^\ast - x^\ast \| ^ 2$ of our DP-Dishuf-AC-based solver (Algorithm 3) and DP-GT-based solver (Algorithm 1) under varying network size $n=10,50,250$. To make a comparison, we consider a differentially private average consensus based solver (DP-AC-based solver), obtained by removing the DiShuf mechanism in Algorithm 3, where Gaussian noises are directly added to privacy-sensitive data to preserve differential privacy \cite{ImprovedGaussianDP}. It can be observed in Figure \ref{fig:n comparison} that our Algorithm 3 (i.e., the DP-Dishuf-AC-based solver) shows the best computation accuracy in different scale networks.

\section{Conclusion}
\label{sec-6}
In this paper, two solvers were proposed to solve the distributed least squares optimization problem with requirements of preserving differential privacy of local cost functions. The first solver was established on the distributed gradient-tracking algorithm,  by perturbing the parameters encoding the sensitive data with truncated Laplacian noises and Gaussian noises. It was shown that this method leads to a limited differential-privacy-protection ability and the resulting computation accuracy linearly relies on the network size. To overcome these issues, the second solver was developed where the differentially private DiShuf-based average consensus algorithm was employed such that each agent obtains an noisy estimate of the global gradient function and thus compute the solution directly. As a result, one could achieve arbitrary differential privacy requirements  and a computation accuracy independent of the network size. By simulations, it was shown that the computation accuracy of the former solver strictly increases as the network size increases, while such trend was not observed by the latter.

\appendix

\subsection{Proof of Theorem 1}
\emph{i).} 
For each agent $i\in\mathrm{V}$, consider the following two mechanisms:
\begin{equation}
    \begin{aligned}
        & \mathcal{M}_i^A=\theta_i^A+\gamma_i \\
        & \mathcal{M}_i^B=B_i+\eta_i.
    \end{aligned}
\end{equation}
By \cite[Theorem 1]{BoundedNoiseDP}, the truncated Laplacian mechanism $\mathcal{M}_i^A$ can be easily verified to be $(\epsilon,\delta)$-differentially private under $\mu$ adjacency, with $(\epsilon,\delta,\mu)$ satisfying \eqref{eq:DP-GT slover delta} and each entry of $\gamma_i$ generated according to \eqref{eq:Truncated Laplacian Distribution}. Similarly, by \cite[Theorem 8]{ImprovedGaussianDP}, the Gaussian mechanism $\mathcal{M}_i^B$ can be found to be $(\epsilon,\delta)$-differentially private under $\mu$ adjacency, with $\eta_i\sim\mathcal{N}(0,\sigma_\eta^2)^m$ and $\sigma_\eta \geq \mu/{\bar\kappa_{\epsilon}(\delta)}$. With these properties in mind, we then turn to take a look at all communication messages that are eavesdropped during computation, and notice that they are in fact deterministic mappings of all mechanisms $\mathcal{M}_i^A,\mathcal{M}_i^B$, $i\in\mathrm{V}$. Thus, by employing the parallel composition \cite[Theorem 4]{DPComposition} and post-processing \cite[Proposition 2.1]{DPpost-processing} property of differential privacy, we conclude that Algorithm 1 preserves $(\epsilon,\delta)$-differential privacy of $\theta_i, i\in\mathrm{V}$.

\emph{ii).} 
It is noted from \cite{GradientTracking} that Algorithm 1 indeed solves the least-squares optimization problem of the from 
\[
\min_{x\in \mathbb{R}^{m}} \bar f(x)=\sum_{i=1}^{n}\frac{1}{2}x^{\top}G_{i}x+H_{i}^{\top}x\,.
\]
More explicitly, by \cite[Theorem 1]{GradientTracking} there exists $\beta^\ast>0$ such that each state $x_i(t)$  converges to the unique solution if it exists, with some linear convergence rate $\alpha_1\in(0,1)$. With this in mind, we let $G=\sum_{i=1}^{n}G_{i}$ and $H=\sum_{i=1}^{n}H_{i}$, which implies $G=A+\Omega_{A}$ and $H=B+\Omega_{B}$.
By noting that each element of $\Omega_A$ is a truncated Laplacian noise according to \eqref{eq:Truncated Laplacian Distribution}, it can be seen that $\|\Omega_{A}\| \leq \|\Omega_{A}\|_{F} < \sqrt{n}m\bar\gamma$. Then, with Assumption \ref{ass-1} there holds  $G=A+\Omega_{A}\geq (1-d)\underline{\lambda}_{A}\mathbf{I}_m>0$, implying that $\bar f$ is strongly convex and admits the unique solution $x(\infty)=-G^{-1}H$. The statement ii) is thus completed.

\emph{iii).} 
We note that
\begin{equation}
    \begin{aligned}
        & \|x(\infty)-x^{\ast}\|^2 \\
        = & \|(A+\Omega_A)^{-1}B + (A+\Omega_A)^{-1}\Omega_B - A^{-1}B\|^2 \\
        \leq & 2\|(A+\Omega_A)^{-1}B - A^{-1}B\|^2 + 2\|(A+\Omega_A)^{-1}\Omega_B\|^2 \\
        \leq & 2\|((A+\Omega_A)^{-1}-A^{-1})B\|^2 + 2\|(A+\Omega_A)^{-1}\Omega_B\|^2 \\
        = & 2\|(A+\Omega_A)^{-1}[I-(A+\Omega_A)A^{-1}]B\|^2 \\
        & + 2\|(A+\Omega_A)^{-1}\Omega_B\|^2 \\
        = & 2\|(A+\Omega_A)^{-1}\Omega_AA^{-1}B\|^2 + 2\|(A+\Omega_A)^{-1}\Omega_B\|^2 \\
        = & 2\|(A+\Omega_A)^{-1}\Omega_Ax^{\ast}\|^2 + 2\|(A+\Omega_A)^{-1}\Omega_B\|^2\\
        \leq & 2\|(A+\Omega_A)^{-1}\|^2(\|\Omega_Ax^{\ast}\|^2 + \|\Omega_B\|^2)\,.
    \end{aligned}
\end{equation}
By recalling that $A+\Omega_{A}\geq (1-d)\underline{\lambda}_{A}\mathbf{I}_m$, we have $\|(A+\Omega_A)^{-1}\|\leq \frac{1}{(1-d)\underline{\lambda}_{A}}$. This yields
\[\begin{array}{rcl}
\|x(\infty)-x^{\ast}\|^2 &\leq& \frac{2}{(1-d)^2\underline{\lambda}_{A}^2}(\|\Omega_Ax^{\ast}\|^2 + \|\Omega_B\|^2) \\
&\leq& \frac{2}{(1-d)^2\underline{\lambda}_{A}^2}(\|x^{\ast}\|^2 \|\Omega_A\|_F^2 + \|\Omega_B\|^2)\,.
\end{array}\]
To this end, we proceed to study the bound of $\mathbb{E}\|\Omega_A\|_F$ and $\mathbb{E}\|\Omega_B\|^2$, and have the following observations.
\begin{itemize}
    \item  Since $\Omega_A=\sum_{i=1}^{n}\mathcal{F}^{-1}(\gamma_{i})$, we can obtain that each element of matrix $\Omega_A$ is the sum of $n$ i.i.d. truncated Laplacian noises with zero mean and $\sigma_\gamma^2$ variance. Thus, $\mathbb{E}\|\Omega_{A}\|^2 \leq \mathbb{E}\|\Omega_{A}\|^2_{F} = nm^2\sigma_\gamma^2$.
    \item  Similarly, each element of $\Omega_{B}$ is the sum of $n$ i.i.d. noises subject to $\mathcal{N}(0,\sigma_\eta^2)$. Then $\mathbb{E} \|\Omega_{B}\|^2 = mn\sigma_\eta^2$.
\end{itemize}
Thus, we have
\[
    \begin{aligned}
        & \mathbb{E}\|x(\infty)-x^{\ast}\|^2 \\
        \leq & \frac{2}{(1-d)^2\underline{\lambda}_{A}^2}(\|x^{\ast}\|^2 \mathbb{E}\|\Omega_A\|_F^2 + \mathbb{E}\|\Omega_B\|^2)\\
        \leq & \frac{2}{(1-d)^2\underline{\lambda}_{A}^2}(nm^2\sigma_\gamma^2 \|x^{\ast}\|^2 + nm\sigma_\eta^2).
    \end{aligned}
\]
The proof is thus completed.

\bibliographystyle{IEEEtran}
\bibliography{ref}
\end{document}